\documentclass[twoside]{amsart}
\usepackage{amsfonts,amsthm}
\usepackage{enumerate}
\theoremstyle{plain}
\newtheorem*{thm A}{Theorem~A}
\newtheorem*{thm B}{Theorem~B}
\newtheorem*{thm C}{Theorem~C}
\newtheorem*{thm D}{Theorem~D}
\newtheorem*{thm E}{Theorem~E}
\newtheorem*{mthm}{Main Theorem}
\newtheorem*{thm 1}{Theorem~1}
\newtheorem*{thm 2}{Theorem~2}
\newtheorem*{pro A}{Proposition~A}
\newtheorem*{pro B}{Proposition~B}
\newtheorem*{lem A}{Lemma~A}
\newtheorem*{lem B}{Lemma~B}
\newtheorem*{lem C}{Lemma~C}
\newtheorem*{lem D}{Lemma~D}

\newtheorem*{proof*}{\it The proof of Corollary}

\newtheorem{theorem}{Theorem}[section]
\newtheorem{corollary}[theorem]{Corollary}
\newtheorem{lemma}[theorem]{Lemma}
\newtheorem{proposition}[theorem]{Proposition}
\newtheorem{remark}[theorem]{Remark}

\theoremstyle{plain}
\newcommand{\be}{\begin{equation}}
\newcommand{\ee}{\end{equation}}
\newcommand{\bea}{\begin{eqnarray}}
\newcommand{\eea}{\end{eqnarray}}
\newcommand{\ba}{\begin{array}}

\newcommand{\ea}{\end{array}}

\newcommand{\bc}{\begin{center}}
\newcommand{\ec}{\end{center}}

\newcommand{\benu}{\begin{enumerate}}
\newcommand{\eenu}{\end{enumerate}}
\newcommand{\bpr}{\begin{proposition}}
\newcommand{\epr}{\end{proposition}}
\newcommand{\ble}{\begin{lemma}}
\newcommand{\ele}{\end{lemma}}
\newcommand{\bco}{\begin{corollary}}
\newcommand{\eco}{\end{corollary}}

\def \Cal{\mathcal}

\def \GBt{G_2({\mathbb C}^{m+2})}

\def \EK{{\eta}({\xi})}
\def \EN{{\eta}_{\nu}}
\def \EtK{{\eta}_2(\xi)}
\def \EsK{{\eta}_3(\xi)}
\def \ENK{{\eta}_{\nu}({\xi})}
\def \ENoK{{\eta}_{{\nu}+1}({\xi})}
\def \ENtK{{\eta}_{{\nu}+2}({\xi})}
\def \EKN{{\eta}({\xi}_{\nu})}
\def \ENSK{{\eta}_{\nu}^2({\xi})}

\def \Et{{\eta}_2}
\def \Es{{\eta}_3}
\def \Kt{{\xi}_2}
\def \Ks{{\xi}_3}

\def \Na{\nabla}
\def \SN{{\sum}_{{\nu}=1}^3}

\def \PN{{\phi}_{\nu}}
\def \Po{{\phi}_1}
\def \Pt{{\phi}_2}
\def \Ps{{\phi}_3}
\def \PoP{{\phi}_1{\phi}}
\def \PtP{{\phi}_2{\phi}}

\def \PtK{{\phi}_2{\xi}}
\def \PsK{{\phi}_3{\xi}}

\def \PKt{{\phi}{\xi}_2}
\def \PKs{{\phi}{\xi}_3}

\def \PNK{{\phi}_{\nu}{\xi}}
\def \PKN{{\phi}{\xi}_{\nu}}
\def \PKNo{{\phi}{\xi}_{{\nu}+1}}
\def \PKNt{{\phi}{\xi}_{{\nu}+2}}

\def \PPN{{\phi}{\phi}_{\nu}}

\def \PPo{{\phi}{\phi}_1}

\def \PPoK{{\phi}{\phi}_1{\xi}}
\def \PPtK{{\phi}{\phi}_2{\xi}}
\def \PPsK{{\phi}{\phi}_3{\xi}}

\def \PPNK{{\phi}{\phi}_{\nu}{\xi}}
\def \PPNP{{\phi}{\phi}_{\nu}{\phi}}
\def \PY{{\phi}AY}

\def \PPNtK{{\phi}{\phi}_{{\nu}+2}{\xi}}
\def \PNo{{\phi}_{{\nu}+1}}
\def \PNt{{\phi}_{{\nu}+2}}
\def \PNP{{\phi}_{\nu}{\phi}}
\def \PNoP{{\phi}_{{\nu}+1}{\phi}}
\def \PNtP{{\phi}_{{\nu}+2}{\phi}}
\def \ENo{{\eta}_{{\nu}+1}}
\def \ENt{{\eta}_{{\nu}+2}}
\def \qNo{q_{{\nu}+1}}
\def \qNt{q_{{\nu}+2}}
\def \qNoK{q_{{\nu}+1}({\xi})}
\def \qNtK{q_{{\nu}+2}({\xi})}

\def \KN{{\xi}_{\nu}}

\def \Ko{{\xi}_1}
\def \Kt{{\xi}_2}
\def \Ks{{\xi}_3}
\def \KNo{{\xi}_{{\nu}+1}}
\def \KNt{{\xi}_{{\nu}+2}}

\def \SI{{\sum}_{i=1}^{4m-1}}
\def \NBo{SU_{2,m-1}/S(U_2{\cdot}U_{m-1})}
\def \NBt{SU_{2,m}/S(U_{2}{\cdot}U_{m})}
\def \CC{{\Bbb C}}
\def \RR{{\Bbb R}}
\def \HH{{\Bbb H}}
\def \RF{\it Reeb flow}

\def \EN{{\eta}_{\nu}}
\def \Na{\nabla}
\def \SN{\sum_{{\nu}=1}^3}

\def \PN{{\phi}_{\nu}}
\def \PPN{{\phi}{\phi}_{\nu}}
\def \PNo{{\phi}_{{\nu}+1}}
\def \PNt{{\phi}_{{\nu}+2}}
\def \PNP{{\phi}_{\nu}{\phi}}
\def \ENo{{\eta}_{{\nu}+1}}
\def \ENt{{\eta}_{{\nu}+2}}
\def \qNo{q_{{\nu}+1}}
\def \qNt{q_{{\nu}+2}}

\def \KN{{\xi}_{\nu}}
\def \KNo{{\xi}_{{\nu}+1}}
\def \KNt{{\xi}_{{\nu}+2}}

\def \SI{{\sum}_{i=1}^{4m-1}}
\def \GBt{G_2({\Bbb C}^{m+2})}
\def \GBo{G_2({\Bbb C}^{m+1})}

\def \SI{{\sum}_{i=1}^{4m-1}}
\def \NBo{SU_{2,m-1}/S(U_2{\cdot}U_{m-1})}
\def \NBt{SU_{2,m}/S(U_{2}{\cdot}U_{m})}
\def \CC{{\Bbb C}}
\def \RR{{\Bbb R}}
\def \HH{{\Bbb H}}
\begin{document}

\title[Hypersurfaces in
complex hyperbolic two-plane Grassmannians]{Real Hypersurfaces in Complex Hyperbolic Two-Plane
Grassmannians with commuting\\ Ricci tensor}
\vspace{0.2in}
\author[Young Jin Suh]{Young Jin Suh}
\address{\newline
Young Jin Suh
\newline Kyungpook National University,
\newline Department of Mathematics,
\newline Taegu 702-701, Korea}
\email{yjsuh@knu.ac.kr}


\thanks{* This work was supported by grants Proj. Nos. NRF-2012-R1A2A2A-01043023
from National Research Foundation of Korea.}

\begin{abstract}
In this paper we first
introduce the full expression of the curvature tensor of a real
hypersurface $M$ in complex hyperbolic two-plane Grassmannians $SU_{2,m}/S(U_2{\cdot}U_m)$, $m{\ge}2$ from the equation of Gauss. Next we derive a new formula
for the Ricci tensor of $M$ in $SU_{2,m}/S(U_2{\cdot}U_m)$. Finally we give a complete classification of Hopf hypersurfaces in
complex hyperbolic two-plane Grassmannians $SU_{2,m}/S(U_2{\cdot}U_m)$ with commuting
Ricci tensor.  Each can be described as a tube over
a totally geodesic $SU_{2,m-1}/S(U_2{\cdot}U_{m-1})$ in
$SU_{2,m}/S(U_2{\cdot}U_m)$ or a horosphere whose center at infinity
is singular.
\end{abstract}

\maketitle

\section*{Introduction}
\setcounter{equation}{0}
\renewcommand{\theequation}{0.\arabic{equation}}
\vspace{0.13in}

In the geometry of real hypersurfaces in complex space forms $M_m(c)$ or in
quaternionic space forms $Q_m(c)$ Kimura \cite{K} and \cite{K2}
 (resp. P\'erez and the author \cite{PS2}) considered real hypersurfaces
in $M_n(c)$ (resp. in $Q_m(c)$) with commuting Ricci tensor, that is, $S{\phi}={\phi}S$, (resp. $S{\phi}_i={\phi}_iS$, $i=1,2,3$)where $S$ and $\phi$
(resp. $S$ and ${\phi}_i$, $i=1,2,3$) denote the Ricci tensor and the structure tensor of real hypersurfaces in $M_m(c)$ (resp. in $Q_m(c)$).
\vskip 6pt
\par
In \cite{K} and \cite{K2}, Kimura has classified that a Hopf hypersurface $M$ in complex projective space $P_m({\Bbb C})$
with commuting Ricci tensor is locally congruent to of type $(A)$, a tube over a totally geodesic $P_k({\Bbb C})$, of type $(B)$, a tube over a complex quadric
$Q_{m-1}$, $\cot^2 2r=m-2$, of type $(C)$, a tube over
$P_1({\Bbb C}){\times}P_{(m-1)/2}({\Bbb C})$, $\cot^2 2r=\frac{1}{m-2}$ and $n$ is odd, of type $(D)$, a tube over a
complex two-plane Grassmannian $G_2({\Bbb C}^5)$, $\cot^22r=\frac{3}{5}$ and $n=9$, of type $(E)$, a tube over a Hermitian symmetric space $SO(10)/U(5)$,
$\cot^22r=\frac{5}{9}$ and $m=15$.

\par
\vskip 6pt
On the other hand, in a quaternionic projective space
${\Bbb Q}P^m$ P\'erez and the author \cite{PS2} have classified real hypersurfaces in
$QP^m$ with commuting Ricci tensor $S{\phi}_i={\phi}_iS$,$i=1,2,3$,
where $S$ (resp. $\phi_i$) denotes the Ricci tensor(resp. the
structure tensor) of $M$ in ${\Bbb Q}P^m$, is locally congruent to
of $A_1, A_2$-type, that is, a tube over ${\Bbb Q}P^k$ with radius
$0<r<\frac{\pi}{2}$, $k{\in}\{0,{\cdots},m-1\}$. The almost contact
structure vector fields $\{\xi_1 ,\xi_2 ,\xi_3 \}$ are defined by
$\xi_i =-J_i N$, $i=1,2,3$, where $J_i$, $i=1,2,3$, denote a
quaternionic K\"ahler structure of ${\Bbb Q}P^m$ and $N$ a unit
normal field of $M$ in ${\Bbb Q}P^m$.  Moreover, P\'erez and Suh \cite{PS} have considered the notion of
${\nabla}_{\xi_i}R=0$, $i=1,2,3$, where $R$ denotes the curvature
tensor of a real hypersurface $M$ in ${\Bbb Q}P^m$, and proved that
$M$ is locally congruent to a tube of radius $\frac{\pi}{4}$ over
${\Bbb Q}P^k$.
\par
\vskip 6pt

Let us denote by $SU_{2,m}$ the set of
$(m+2){\times}(m+2)$-indefinite special unitary matrices
 and $U_m$ the set of $m{\times}m$-unitary matrices. Then
the Riemannian symmetric space $SU_{2,m}/S(U_2U_m)$, $m \geq 2$,
which consists of complex two-dimensional subspaces in indefinite
complex Euclidean space ${\Bbb C}_2^{m+2}$, has a remarkable
feature that it is a Hermitian symmetric space as well as a
quaternionic K\"{a}hler symmetric space. In fact, among all
Riemannian symmetric spaces of noncompact type the symmetric spaces
$SU_{2,m}/S(U_2U_m)$, $m \geq 2$, are the only ones which are
Hermitian symmetric and quaternionic K\"{a}hler symmetric.
\par
\vskip 6pt
The existence of these two structures leads to a number of interesting
geometric problems on $SU_{2,m}/S(U_2U_m)$, one of which we are
going to study in this article. To describe this problem, we denote
by $J$ the K\"{a}hler structure and by ${\frak J}$ the quaternionic
K\"{a}hler structure on $SU_{2,m}/S(U_2U_m)$. Let $M$ be a connected
hypersurface in $SU_{2,m}/S(U_2U_m)$ and denote by $N$ a unit normal to $M$.
Then a structure vector field $\xi$ defined by ${\xi}=-JN$ is said to be a Reeb vector field.
\par
\vskip 6pt

Next, we consider the standard embedding of $SU_{2,m-1}$ in
$SU_{2,m}$. Then the orbit $SU_{2,m-1} \cdot o$ of $SU_{2,m-1}$
through $o$ is the Riemannian symmetric space \hfill
\newline
$SU_{2,{m-1}}/S(U_2U_{m-1})$ embedded in $SU_{2,m}/S(U_2U_m)$ as a
totally geodesic submanifold. Every tube around
$SU_{2,m-1}/S(U_2U_{m-1})$ in $SU_{2,m}/S(U_2U_m)$ has the property
that both maximal complex subbundle ${\Cal C}$ and quaternionic subbundle ${\Cal Q}$ are invariant under the shape
operator.
\par
\vskip 6pt
Finally, let $m$ be even, say $m = 2n$, and consider the standard
embedding of $Sp_{1,n}$ in $SU_{2,2n}$. Then the orbit $Sp_{1,n}
\cdot o$ of $Sp_{1,n}$ through $o$ is the quaternionic hyperbolic
space ${\Bbb H}H^n$ embedded in $SU_{2,2n}/S(U_2U_{2n})$ as a
totally geodesic submanifold. Any tube around ${\Bbb H}H^n$ in
$SU_{2,2n}/S(U_2U_{2n})$ has the property that both ${\Cal C}$ and
${\Cal Q}$ are invariant under the shape operator.
\par
\vskip 6pt
As a converse of the statements mentioned above, we assert that
with one possible exceptional case there are no other such real
hypersurfaces. Related to such a result, we introduce another theorem due to Berndt and Suh \cite{BS3} as follows:
\vskip 6pt
\par
\begin{thm A}\label{Theorem A}
Let $M$ be a connected hypersurface in
$SU_{2,m}/S(U_2U_m)$, $m \geq 2$. Then the maximal complex subbundle
${\Cal C}$ of $TM$ and the maximal quaternionic subbundle ${\Cal Q}$
of $TM$ are both invariant under the shape operator of $M$ if and
only if $M$ is congruent to an open part of one of the following
hypersurfaces:
\par
$(A)$\quad a tube around a totally geodesic
$SU_{2,m-1}/S(U_2U_{m-1})$ in $SU_{2,m}/S(U_2U_m)$;
\par
$(B)$\quad a tube around a totally geodesic ${\Bbb H}H^n$ in
$SU_{2,2n}/S(U_2U_{2n})$, $m = 2n$;
\par
$(C)$\quad a horosphere in $SU_{2,m}/S(U_2U_m)$
 whose center at infinity is singular;
\par
or the following exceptional case holds:
\par
$(D)$\quad The normal bundle $\nu M$ of $M$ consists of singular
tangent vectors of type $JX \perp {\frak J}X$. Moreover, $M$ has at
least four distinct principal curvatures, three of which are given
by
\begin{equation*}
\alpha = \sqrt{2}\ ,\ \gamma = 0\ ,\ \lambda = \frac{1}{\sqrt{2}}
\end{equation*}
with corresponding principal curvature spaces
\begin{equation*}
T_\alpha = TM \ominus ({\Cal C} \cap {\Cal Q})\ ,\ T_\gamma = J(TM
\ominus {\Cal Q})\ ,\ T_\lambda \subset  {\Cal C} \cap {\Cal Q} \cap
J{\Cal Q}.
\end{equation*}
If $\mu$ is another (possibly nonconstant) principal curvature
function, then we have $T_\mu \subset {\Cal C} \cap {\Cal Q} \cap
J{\Cal Q}$, $JT_\mu \subset T_\lambda$ and ${\frak J}T_\mu \subset
T_\lambda $.
\end{thm A}
\par
\vskip 6pt
In Theorem A the maximal complex subbundle $\Cal C$ of $TM$ is
invariant under the shape operator if and only if the Reeb vector field
$\xi$ becomes a principal vector field for the shape operator $A$ of
$M$ in $SU_{2,m}/S(U_2U_m)$. In this case the Reeb vector field
$\xi$ is said to be a Hopf vector field.
The flow generated by the integral curves of the structure
vector field $\xi$ for Hopf hypersurfaces in $G_2({\Bbb C}^{m+2})$
is said to be a {\it geodesic Reeb flow}.
\par
\vskip 6pt
 The classification of all real hypersurfaces in complex
projective space $\CC P^m$ with isometric {\RF} has been obtained by
Okumura \cite{O}. The corresponding classification in complex
hyperbolic space $\CC H^m$ is due to Montiel and Romero \cite{MR}
and in quaternionic projective space $\HH P^m$ due to Martinez and
P\'erez \cite{MP} respectively.
\par
\vskip 6pt

\vskip 6pt
Now let us introduce a classification theorem due to Suh \cite{S5} for all real
hypersurfaces with isometric {\RF} in complex hyperbolic two-plane
Grassmann manifold $\NBt$ as follows: \hfill
\newline
\vskip 6pt
\par
\begin{thm B}\label{Theorem B}
Let $M$ be a connected orientable real hypersurface in the complex
hyperbolic two-plane Grassmannian $\NBt$, $m \geq 3$. Then the Reeb
flow on $M$ is isometric if and only if $M$ is an open part of a
tube around some totally geodesic $\NBo$ in $\NBt$ or a horosphere
whose center at infinity with $JX{\in}{\frak J}X$ is singular.
\end{thm B}
\vskip 6pt

In the proof of Theorem A we proved that the
1-dimensional distribution $[\xi ]$ is contained in either the
3-dimensional distribution ${\mathcal Q}^{\bot}$ or in the orthogonal
complement ${\mathcal Q}$ such that $T_xM={\mathcal Q}{\oplus}{\mathcal Q}^{\bot}$.
The case $(A)$ in Theorem A is just the case that the 1-
dimensional distribution $[\xi ]$ belongs to the distribution
${\mathcal Q}$. Of course, it is not difficult to check that the
Ricci tensor $S$ of type $(A)$ or of type $(C)$ with $JX{\in}{\frak J}X$ in Theorem A commutes with the structure tensor, that is $S{\phi}={\phi}S$.  Then it must be a natural question to ask whether real
hypersurfaces in $\NBt$ with commuting Ricci tensor
can exist or not . \vskip 6pt
\par

\vskip 6pt
In this paper we consider such a converse problem and want to give a complete classification of real hypersurfaces in $G_2({\Bbb C}^{m+2})$ satisfying $S{\phi}={\phi}S$ as follows:
\vskip 6pt
\par
\begin{mthm}\label{Main Theorem}
Let $M$ be a Hopf hypersurface in $\NBt$ with commuting Ricci tensor, $m{\ge}3$. Then $M$ is locally congruent to an open part of a
tube around some totally geodesic $\NBo$ in $\NBt$ or a horosphere
whose center at infinity with $JX{\in}{\frak J}X$ is singular.
\end{mthm}
\vskip 6pt
\par

\vskip 6pt
 A remarkable consequence of
our Main Theorem is that a connected complete real hypersurface in
$\NBt$, $m \geq 3$ with commuting Ricci tensor is homogeneous and has an isometric Reeb flow. This was
also true in complex two-plane Grassmannians $G_2({\Bbb C}^{m+2})$,
which could be identified with symmetric space of compact type
$SU_{m+2}/S(U_2{\cdot}U_m)$, as follows from the classification. It
would be interesting to understand the actual reason for it (See
\cite{BS1}, \cite{S1} and \cite{S2}).
\par
\vskip 6pt

This paper is organized as follows. In Section 1 we summarize some
basic facts about the Riemannian geometry of $\NBt$. In Section 2 we
obtain some basic geometric equations for real hypersurfaces in $\NBt$.
In Section~\ref{section 3} we study real hypersurfaces in $\NBt$ with Ricci commuting for ${\xi}{\in}{\mathcal Q}$ , and in Section~\ref{section 4} those with Ricci commuting for ${\xi}{\in}{\mathcal Q}^{\bot}$. Finally in
Section~\ref{section 5} we use these results to derive our classification.

\vskip 8pt

\section  {The complex hyperbolic two-plane Grassmannian $\NBt$} \setcounter{equation}{0}
\renewcommand{\theequation}{1.\arabic{equation}}
\vspace{0.13in}

In this section we summarize basic material about complex hyperbolic
Grassmann manifolds $\NBt$, for details we refer to
\cite{BCO}, \cite{BS3},\cite{S2} and \cite{S5}.
\par
\par
\vskip 6pt The Riemannian symmetric space
$SU_{2,m}/S(U_2{\cdot}U_m)$, which consists of all complex
two-dimensional linear subspaces in indefinite complex Euclidean
space $\CC_2^{m+2}$, becomes a connected, simply connected,
irreducible Riemannian symmetric space of noncompact type and with
rank two. Let $G = SU_{2,m}$ and $K = S(U_2{\cdot}U_m)$, and denote
by ${\frak g}$ and ${\frak k}$ the corresponding Lie algebra of the
Lie group $G$ and $K$ respectively. Let $B$ be the Killing form of
${\frak g}$ and denote by ${\frak p}$ the orthogonal complement of
${\frak k}$ in ${\frak g}$ with respect to $B$. The resulting
decomposition ${\frak g} = {\frak k} \oplus {\frak p}$ is a Cartan
decomposition of ${\frak g}$. The Cartan involution $\theta \in
{\text Aut}({\frak g})$ on ${\frak s}{\frak u}_{2,m}$ is given by
$\theta(A) = I_{2,m} A I_{2,m}$, where
\begin{equation*}
\begin{split}
I_{2,m} = \begin{pmatrix}
-I_{2} & 0_{2,m} \\
0_{m,2} & I_{m} \end{pmatrix}
\end{split}
\end{equation*}
$I_2$ and $I_m$ denotes the identity $(2 \times 2)$-matrix and $(m
\times m)$-matrix respectively. Then $< X , Y > = -B(X,\theta Y)$
becomes a positive definite ${\text Ad}(K)$-invariant inner product
on ${\frak g}$. Its restriction to ${\frak p}$ induces a metric $g$
on $SU_{2,m}/S(U_2{\cdot}U_m)$, which is also known as the Killing
metric on $SU_{2,m}/S(U_2{\cdot}U_m)$. Throughout this paper we
consider $SU_{2,m}/S(U_2{\cdot}U_m)$ together with this particular
Riemannian metric $g$.
\par
\vskip 6pt The Lie algebra ${\frak k}$ decomposes orthogonally into
${\frak k}  = {\frak s}{\frak u}_2 \oplus {\frak s}{\frak u}_m
\oplus {\frak u}_1$, where ${\frak u}_1$ is the one-dimensional
center of ${\frak k}$. The adjoint action of ${\frak s}{\frak u}_2$
on ${\frak p}$ induces the quaternionic K\"{a}hler structure ${\frak
J}$ on $SU_{2,m}/S(U_2{\cdot}U_m)$, and the adjoint action of
\begin{equation*}
\begin{split}
Z = \begin{pmatrix}
 \frac{mi}{m+2}I_2 & 0_{2,m} \\
 0_{m,2} & \frac{-2i}{m+2}I_m
 \end{pmatrix}
 \in {\frak u}_1
\end{split}
\end{equation*}
induces the K\"{a}hler structure $J$ on $SU_{2,m}/S(U_2{\cdot}U_m)$.
\vskip 6pt
\par

By construction, $J$ commutes with each almost Hermitian structure
$J_{\nu}$ in ${\frak J}$ for ${\nu}=1,2,3$. Recall that a canonical
local basis $J_1,J_2,J_3$ of a quaternionic K\"{a}hler structure
${\frak J}$ consists of three almost Hermitian structures
$J_1,J_2,J_3$ in ${\frak J}$ such that $J_\nu J_{\nu+1} = J_{\nu +
2} = - J_{\nu+1} J_\nu$, where the index $\nu$ is to be taken modulo
$3$. The tensor field $JJ_\nu$, which is locally defined on
$SU_{2,m}/S(U_2{\cdot}U_m)$, is selfadjoint and satisfies
$(JJ_\nu)^2 = I$ and ${\text tr}(JJ_\nu) = 0$, where $I$ is the
identity transformation. For a nonzero tangent vector $X$ we define
${\Bbb R}X = \{\lambda X \vert \lambda \in {\Bbb R}\}$, ${\Bbb C}X =
{\Bbb R}X \oplus {\Bbb R}JX$, and ${\Bbb H}X = {\Bbb R}X \oplus
{\frak J}X$.
\par
\vskip 6pt
We identify the tangent space $T_oSU_{2,m}/S(U_2{\cdot}U_m)$ of
$SU_{2,m}/S(U_2{\cdot}U_m)$ at $o$ with ${\frak p}$ in the usual way. Let
${\frak a}$ be a maximal abelian subspace of ${\frak p}$. Since
$SU_{2,m}/S(U_2{\cdot}U_m)$ has rank two, the dimension of any such
subspace is two. Every nonzero tangent vector $X \in
T_oSU_{2,m}/S(U_2{\cdot}U_m) \cong {\frak p}$ is contained in some maximal
abelian subspace of ${\frak p}$. Generically this subspace is
uniquely determined by $X$, in which case $X$ is called regular.
\par
\vskip 6pt
If there exists more than one maximal abelian subspaces of ${\frak p}$
containing $X$, then $X$ is called singular. There is a simple and
useful characterization of the singular tangent vectors: A nonzero
tangent vector $X \in {\frak p}$ is singular if and only if $JX \in
{\frak J}X$ or $JX \perp {\frak J}X$.
\par
\vskip 6pt
 Up to scaling there exists a unique $S(U(2) \cdot
U(m))$-invariant Riemannian metric $g$ on $\NBt$. Equipped with this
metric $\NBt$ is a Riemannian symmetric space of rank two which is
both K\"ahler and quaternionic K\"ahler.
\par
\vskip 6pt
For computational reasons
we normalize $g$ such that the minimal sectional curvature of
$({\NBt},g)$ is $-4$. The sectional curvature $K$ of the noncompact
symmetric space $SU_{2,m}/S(U_2{\cdot}U_m)$ equipped with the
Killing metric $g$ is bounded by $-4{\leq}K{\leq}0$. The sectional
curvature $-4$ is obtained for all $2$-planes ${\Bbb C}X$ when $X$
is a non-zero vector with $JX{\in}{\frak J}X$.
\par
\vskip 6pt

When $m=1$,  $G_2^{*}(\CC^3)=SU_{1,2}/S(U_1{\cdot}U_2)$ is isometric
to the two-dimensional complex hyperbolic space $\CC H^2$ with
constant holomorphic sectional curvature $-4$.
\par
\vskip 6pt
When $m=2$, we note that the isomorphism $SO(4,2)\simeq SU(2,2)$ yields an isometry between $G_2^{*}(\CC ^4)=SU_{2,2}/S(U_2{\cdot}U_2)$ and
the indefinite real Grassmann manifold $G_2^{*}(\RR_2^6)$ of oriented
two-dimensional linear subspaces of an indefinite Euclidean space $\RR_2^6$.
For this reason we
assume $m \geq 3$ from now on, although many of the subsequent
results also hold for $m = 1,2$.
\par
\vskip 6pt

The Riemannian curvature tensor $\bar{R}$ of $\NBt$ is locally given
by
\begin{equation}\label{1.1}
\begin{split}
\bar{R}(X,Y)Z = &- \frac{1}{2}\Big[ g(Y,Z)X - g(X,Z)Y  +  g(JY,Z)JX\\
& - g(JX,Z)JY- 2g(JX,Y)JZ \\
& + \SN \{g(J_\nu Y,Z)J_\nu X - g(J_\nu X,Z)J_{\nu}Y\\
& - 2g(J_\nu X,Y)J_\nu Z\} \\
& + \SN \{g(J_\nu JY,Z)J_\nu JX - g(J_\nu JX,Z)J_{\nu}JY\}\Big ] ,
\end{split}
\end{equation}
where $J_1,J_2,J_3$ is any canonical local basis of ${\frak J}$.
\par
\vskip 6pt

Recall that a maximal flat in a Riemannian symmetric space $\bar{M}$
is a connected complete flat totally geodesic submanifold of
maximal dimension.
A non-zero tangent vector $X$ of $\bar{M}$ is singular if
$X$ is tangent to more than one maximal flat in $\bar{M}$, otherwise $X$
is regular.
The singular tangent vectors of $\NBt$ are precisely
the eigenvectors
and the asymptotic vectors of the self-adjoint endomorphisms $JJ_1$, where
$J_1$ is any almost Hermitian structure in $\frak J$.
In other words, a tangent vector $X$ to $\NBt$ is singular if
and only if $JX \in \frak J X$ or $JX \bot \frak J X$.
\par
\vskip 6pt

Now we want to focus on a singular vector $X$ of type $JX{\in}{\frak
J}X$. In this paper, we will have to compute explicitly Jacobi
vector fields along geodesics whose tangent vectors are all singular
of type $JX \in \frak J X$. For this we need the eigenvalues and
eigenspaces of the Jacobi operator $\bar{R}_X := \bar{R}(.,X)X$. Let
$X$ be a singular unit vector tangent to $\NBt$ of type $JX \in
\frak J X$. Then there exists an almost Hermitian structure $J_1$ in
$\frak J$ such that $JX = J_1X$ and the eigenvalues, eigenspaces and
multiplicities of $\bar{R}_X$ are respectively given by
\vskip 12pt
\medskip
\begin{center}
\begin{tabular}{|c|c|c|}
\hline
\mbox{principal curvature} & \mbox{eigenspace}  & \mbox{multiplicity}  \\
\hline
$0$ & $\RR X \oplus \{Y \vert Y \perp \HH X,\ JY = -J_1Y\}$ & $2m-1$ \\
$-1$ & $\HH X \ominus \CC X \oplus \{Y \vert Y \perp \HH X, \ JY = J_1Y\}$ & $2m$ \\
$-4$ & $\RR JX$ & $1$\\
\hline
\end{tabular}
\end{center}
\medskip
\vskip 12pt
where ${\RR}X$, ${\CC}X$ and ${\HH}X$ denotes the real, complex and quaternionic span of $X$, respectively, and $\CC ^\perp X$ the orthogonal complement of $\CC X$ in $\HH X$.
\par
\par
\vskip 6pt The maximal totally geodesic submanifolds in complex
hyperbolic two-plane Grassmannian $\NBt$ are $\NBo$, $\CC H^m$, $\CC H^k \times \CC H^{m-k}$ ($1 \leq k
\leq [m/2]$), $G_2^{*}(\RR^{m+2})$ and $\HH H^n$ (if $m = 2n$).
The first three are complex
submanifolds and the other two are real submanifolds with respect to the
K\"ahler structure $J$.
The tangent spaces of
the totally geodesic $\CC H^m$ are precisely
the maximal linear subspaces of the
form $\{X \vert JX = J_1X\}$ with some fixed
almost Hermitian structure $J_1 \in \frak J$.

\vskip 8pt

\section  {Real hypersurfaces in $\NBt$} \setcounter{equation}{0}
\renewcommand{\theequation}{2.\arabic{equation}}
\vspace{0.13in}

Let $M$ be a real hypersurface in $\NBt$, that is, a hypersurface in
$\NBt$ with real codimension one. The induced Riemannian metric on
$M$ will also be denoted by $g$, and $\nabla$ denotes the Levi
Civita covariant derivative of $(M,g)$. We denote by $\mathcal C$
and $\mathcal Q$ the maximal complex and quaternionic subbundle of
the tangent bundle $TM$ of $M$, respectively. Now let us put

\begin{equation}\label{2.1}
JX={\phi}X+{\eta}(X)N,\quad J_{\nu}X={\phi}_{\nu}X+{\eta}_{\nu}(X)N
\end{equation}
for any tangent vector field $X$ of a real hypersurface $M$ in
$\NBt$, where ${\phi}X$ denotes the tangential component of $JX$ and
$N$ a unit normal vector field of $M$ in $\NBt$.
\par
\vskip 6pt
From the K\"{a}hler structure $J$ of $\NBt$ there exists
an almost contact metric structure $(\phi,\xi,\eta,g)$ induced on
$M$ in such a way that

\begin{equation}\label{2.2}
{\phi}^2X=-X+{\eta}(X){\xi},\ {\eta}({\xi})=1,\ {\phi}{\xi}=0, \quad
\text{and}\quad {\eta}(X)=g(X,{\xi})
\end{equation}
for any vector field $X$ on $M$ and ${\xi}=-JN$.
\par
\vskip 6pt

 If $M$ is orientable, then the vector field $\xi$ is
globally defined and said to be the induced {\it Reeb vector field}
on $M$. Furthermore, let $J_1,J_2,J_3$ be a canonical local basis of
$\frak J$. Then each $J_\nu$ induces a local almost contact metric
structure $(\phi_\nu,\xi_\nu,\eta_\nu,g)$, ${\nu}=1,2,3$, on $M$.
Locally, $\mathcal C$ is the orthogonal complement in $TM$ of the
real span of $\xi$, and $\mathcal Q$ the orthogonal complement in
$TM$ of the real span of $\{\xi_1,\xi_2,\xi_3\}$.
\par
\vskip 6pt

Furthermore, let $\{J_1,J_2,J_3\}$ be a canonical local basis of
${\frak J}$. Then the quaternionic K\"{a}hler structure $J_\nu$ of
$\NBt$, together with the condition
\begin{equation*}
J_{\nu}J_{\nu+1} = J_{\nu+2} = -J_{\nu+1}J_{\nu}
\end{equation*}
in section $1$, induces an almost contact metric 3-structure
$(\phi_{\nu}, \xi_{\nu}, \eta_{\nu}, g)$ on $M$ as follows:
\begin{equation}\label{2.3}
\begin{split}
&{\phi}_{\nu}^2X=-X+{\eta}_{\nu}(X)({\xi}_{\nu}),\ {\phi}_{\nu}{\xi}_{\nu}=0,\ {\eta}_{\nu}({\xi}_{\nu})=1\\
&{\phi}_{\nu +1}{\xi}_{\nu}=-{\xi}_{{\nu}+2},\quad {\phi}_{\nu
}{\xi}_{{\nu}+1}={\xi}_{{\nu}+2},\\
&{\phi}_{\nu}{\phi}_{{\nu}+1}X
= {\phi}_{{\nu}+2}X+{\eta}_{{\nu}+1}(X){\xi}_{\nu},\\
&{\phi}_{{\nu}+1}{\phi}_{\nu}X=-{\phi}_{{\nu}+2}X+{\eta}_{\nu}(X)
{\xi}_{{\nu}+1}
\end{split}
\end{equation}
for any vector field $X$ tangent to $M$. The tangential and normal
component of the commuting identity $JJ_\nu X = J_\nu JX$ give
\begin{equation}\label{2.4}
\phi\phi_\nu X - \phi_\nu \phi X = \eta_\nu(X)\xi - \eta(X)\xi_\nu \
\text{and}\quad \eta_\nu(\phi X) = \eta(\phi_\nu X).
\end{equation}
The last equation implies $\phi_\nu \xi = \phi \xi_\nu$. The
tangential and normal component of $J_\nu J_{\nu+1}X = J_{\nu+2}X =
-J_{\nu+1}J_\nu X$ give
\begin{equation}\label{2.5}
\phi_\nu\phi_{\nu+1}X - \eta_{\nu+1}(X)\xi_\nu = \phi_{\nu+2}X = -
\phi_{\nu+1}\phi_\nu X + \eta_\nu(X) \xi_{\nu+1}
\end{equation}
and
\begin{equation}\label{2.6}
\eta_\nu(\phi_{\nu+1} X) = \eta_{\nu+2}(X) = - \eta_{\nu+1}(\phi_\nu
X).
\end{equation}

Putting $X = \xi_\nu$ and $X = \xi_{\nu+1}$ into the first of these
two equations yields $\phi_{\nu+2}\xi_\nu = \xi_{\nu+1}$ and
$\phi_{\nu+2}\xi_{\nu+1} = - \xi_\nu$ respectively. Using the Gauss
and Weingarten formulas, the tangential and normal component of the
K\"ahler condition $(\bar{\nabla}_XJ)Y = 0$ give $(\nabla_X\phi)Y =
\eta(Y)AX - g(AX,Y)\xi$ and $(\nabla_X\eta)Y = g(\phi AX,Y)$. The
last equation implies $\nabla_X\xi = \phi AX$. Finally, using the
explicit expression for the Riemannian curvature tensor $\bar{R}$ of
$\NBt$ in \cite{BS3} the Codazzi equation takes the form
\begin{equation}\label{2.7}
\begin{split}
(\nabla_XA)Y &- (\nabla_YA)X
 = -\frac{1}{2}\Big[\eta(X)\phi Y - \eta(Y)\phi X - 2g(\phi X,Y)\xi \\
& \qquad + \SN \big\{\eta_\nu(X)\phi_\nu Y - \eta_\nu(Y)\phi_\nu
 X - 2g(\phi_\nu X,Y)\xi_\nu\big\} \\
& \qquad + \SN \big\{\eta_\nu(\phi X)\phi_\nu\phi Y
- \eta_\nu(\phi Y)\phi_\nu\phi X\big\} \\
& \qquad + \SN \big\{\eta(X)\eta_\nu(\phi Y)
- \eta(Y)\eta_\nu(\phi X)\big\}\xi_\nu \Big] .
\end{split}
\end{equation}
for any vector fields $X$ and $Y$ on $M$. Moreover, by the expression of the curvature tensor \eqref{1.1}, we have the equation of Gauss
as follows:
\begin{equation}\label{2.8}
\begin{split}
R(X,Y)Z  =& -\frac{1}{2}\Bigg[g(Y,Z)X - g(X,Z)Y\\
&+ \ g({\phi}Y,Z){\phi}X - g({\phi}X,Z){\phi}Y - 2g({\phi}X,Y){\phi}Z \\
&+ \SN \left\{g(\PN Y,Z)\PN X - g(\PN X,Z)\PN Y - 2g(\PN X,Y)\PN Z\right\} \\
&+ \SN \left\{g(\PN{\phi}Y,Z)\PN {\phi}X - g(\PN{\phi}X,Z){\PN}{\phi}Y\right\}\\
&- \SN \left\{{\eta}(Y){\EN}(Z){\PN}{\phi}X-{\eta}(X){\EN}(Z){\PN}{\phi}Y\right\}\\
&- \SN \left\{{\eta}(X)g({\PN}{\phi}Y,Z)-{\eta}(Y)g({\PN}{\phi}X,Z)\right\}{\KN}\Bigg]\\
&+g(AY,Z)AX-g(AX,Z)AY
\end{split}
\end{equation}
for any vector fields $X,Y,Z$ and $W$ on $M$.  Hereafter, unless otherwise stated, we want to use these basic
equations mentioned above frequently without referring to them
explicitly.

\bigskip

\section{Some preliminaries in $\NBt$}\label{section 3}
\setcounter{equation}{0}
\renewcommand{\theequation}{3.\arabic{equation}}
\vspace{0.13in}
\par
\vskip 6pt
In this section we can introduce some preliminaries in $\NBt$ corresponding to the formulas given in \cite{S2} from the negative curvature tensor \eqref{2.8}.
Now let us contract $Y$ and $Z$ in the equation of Gauss \eqref{2.8} in section 2. Then the curvature tensor for a real hypersurface $M$ in $\NBt$ gives a Ricci tensor defined by
\begin{equation}\label{3.1}
\begin{split}
SX=&{\SI}R(X,e_i)e_i\\
=&-\frac{1}{2}\Big[(4m+10)X-3{\eta}(X){\xi}-3{\SN}{\EN}(X){\KN}\\
&+ {\SN}\{(\text{Tr}{\PN}{\phi}){\PN}{\phi}X-({\PN}{\phi})^2X\}\\
&-{\SN}\{{\EN}({\xi}){\PN}{\phi}X-{\eta}(X){\PN}{\phi}{\KN}\}\\
&-{\SN}\{(\text{Tr}\ {\PN}{\phi}){\eta}(X)-{\eta}({\PN}{\phi}X)\}{\KN}\Big]+hAX-A^2X ,
\end{split}
\end{equation}
where $h$ denotes the trace of the shape operator $A$ of $M$ in $\NBt$. From the formula $JJ_{\nu}=J_{\nu}J$, $\text{Tr}\ JJ_{\nu}=0$, ${\nu}=1,2,3$ we calculate the following for any basis $\{e_1,{\cdots},e_{4m-1}, N\}$ of the
tangent space of $\NBt$
\begin{equation}\label{3.2}
\begin{split}
0=& \text{Tr}\ JJ_{\nu}\\
=&{\sum}_{k=1}^{4m-1}g(JJ_{\nu}e_k,e_k)+g(JJ_{\nu}N,N)\\
=& \text{Tr}\ {\phi}{\PN} - {\EN}({\xi}) - g(J_{\nu}N,JN)\\
=& \text{Tr}\ {\phi}{\phi}_{\nu} - 2{\EN}({\xi})
\end{split}
\end{equation}
and
\begin{equation}\label{3.3}
\begin{split}
({\PN}{\phi})^2X=&{\PN}{\phi}({\phi}{\PN}X-{\EN}(X){\xi}+{\eta}(X){\KN})\\
=&{\PN}(-{\PN}X+{\eta}({\PN}X){\xi})+{\eta}(X){\PN}^2{\xi}\\
=&X-{\EN}(X){\KN}+{\eta}({\PN}X){\PN}{\xi}\\
&+{\eta}(X)\{-{\xi}+{\EN}({\xi}){\xi}\}.
\end{split}
\end{equation}
Substituting (3.2) and (3.3) into (3.1), we have
\begin{equation}\label{3.4}
\begin{split}
SX=&-\frac{1}{2}\Big[(4m+7)X-3{\eta}(X){\xi}-3{\SN}{\EN}(X){\KN}\\
&+{\SN}\{{\EN}({\xi}){\PN}{\phi}X-{\eta}({\PN}X){\PN}{\xi}-{\eta}(X){\EN}({\xi}){\xi}_{\nu}\}\Big]\\
&+hAX-A^2X .
\end{split}
\end{equation}

Now the covariant derivative of (3.4) becomes
\begin{equation*}
\begin{split}
({\nabla}_YS)X=&\frac{3}{2}(({\nabla}_Y{\eta})X){\xi}+\frac{3}{2}{\eta}(X){\nabla}_Y{\xi}\\
&+\frac{3}{2}{\SN}({\nabla}_Y{\EN})(X){\KN}+\frac{3}{2}{\SN}{\EN}(X){\nabla}_Y{\KN}\\
&-\frac{1}{2}{\SN}\Big\{Y({\EN}({\xi})){\PN}{\phi}X+{\EN}({\xi})({\nabla}_Y{\PN}){\phi}X\\
&+{\EN}({\xi}){\PN}({\nabla}_Y{\phi})X -({\nabla}_Y{\eta})({\PN}X){\PN}{\xi}\\
&-{\eta}(({\nabla}_Y{\PN})X){\PN}{\xi}-{\eta}({\PN}X){\Na}_Y({\PN}{\xi})\\
&-({\Na}_Y{\eta})(X){\EN}({\xi}){\xi}_{\nu}-{\eta}(X){\Na}_Y({\EN}({\xi})){\KN}\\
&-{\eta}(X){\EN}({\xi}){\Na}_Y{\KN}\Big\}\\
&+(Yh)AX+h({\Na}_YA)X-({\Na}_YA^2)X
\end{split}
\end{equation*}
for any vector fields $X$ and $Y$ tangent to $M$ in $\NBt$.
Then from the above formula, together with the formulas in section 2, we have
\begin{equation}\label{3.5}
\begin{split}
({\Na}_YS)X=&\frac{3}{2}g({\phi}AY,X){\xi}+\frac{3}{2}{\eta}(X){\phi}AY\\
&+\frac{3}{2}{\SN}\{{\qNt}(Y){\ENo}(X)-{\qNo}(Y){\ENt}(X)\\
&+g({\PN}AY,X)\}{\KN}\\
&+\frac{3}{2}{\SN}{\EN}(X)\{{\qNt}(Y){\KNo}-{\qNo}(Y){\KNt}+{\PN}AY\}\\
&-\frac{1}{2}{\SN}\Big[Y({\EN}({\xi})){\PNP}X+{\EN}({\xi})\{-{\qNo}(Y){\PNtP}X\\
&+{\qNt}(Y){\PNoP}X + {\EN}({\phi}X)AY-g(AY,{\phi}X){\KN}\}\\
&+{\EN}({\xi})\{{\eta}(X){\PN}AY-g(AY,X){\PN}{\xi}\}-g({\phi}AY,{\PN}X){\PN}{\xi}\\
&+\{{\qNo}(Y){\eta}({\PNt}X)-{\qNt}(Y){\eta}({\PNo}X)-{\EN}(X){\eta}(AY)\\
&+{\eta}({\KN})g(AY,X)\}{\PN}{\xi}\\
&-{\eta}({\PN}X)\{{\qNt}(Y){\PNo}{\xi}-{\qNo}(Y){\PNt}{\xi}+{\PNP}AY\\
&-{\eta}(AY){\KN}+{\eta}({\KN})AY\}\\
&-g({\phi}AY,X){\EN}({\xi}){\KN}-{\eta}(X)Y({\EN}({\xi})){\KN} - {\eta}(X){\EN}({\xi}){\Na}_Y{\KN}\Big]\\
&+(Yh)AX+h({\Na}_YA)X-({\Na}_YA^2)X .
\end{split}
\end{equation}
Now let us take a covariant derivative of $S{\phi}={\phi}S$. Then it gives that
\begin{equation}\label{3.6}
({\nabla}_YS){\phi}X+S({\nabla}_Y{\phi})X=({\nabla}_Y{\phi})SX+{\phi}({\nabla}_YS)X.
\end{equation}
Then the first term of \eqref{3.6} becomes
\begin{equation}\label{3.7}
\begin{split}
({\nabla}_Y&S){\phi}X= \frac{3}{2}g({\phi}AY,{\phi}X){\xi}\\
&+\frac{3}{2}{\SN}\{{\qNt}(Y){\ENo}({\phi}X)-{\qNo}(Y){\ENt}({\phi}X)+g({\PN}AY,{\phi}X)\}{\KN}\\
&+\frac{3}{2}{\SN}{\EN}({\phi}X)\{{\qNt}(Y){\KNo}-{\qNo}(Y){\KNt}+{\PN}A{\phi}X\}\\
&-\frac{1}{2}{\SN}\Big[Y({\EN}({\xi})){\PN}{\phi}^2X + {\EN}({\xi})\{-{\qNo}(Y){\PNt}{\phi}^2X\\
&+{\qNt}(Y){\PNo}{\phi}^2X+{\EN}({\phi}^2X)AY-g(AY,{\phi}^2X){\KN}\}\\
&-{\EN}({\xi})g(AY,{\phi}X){\PNK}-g({\phi}AY,{\PNP}X){\PNK}\\
&+\{{\qNo}(Y){\eta}({\PNt}{\phi}X)-{\qNt}(Y){\eta}({\PNo}{\phi}X)-{\EN}({\phi}X){\eta}(AY)\\
&+{\EKN}g(AY,{\phi}X)\}{\PNK}\\
&-{\eta}({\PNP}X)\{{\qNt}(Y){\PNo}{\xi}-{\qNo}(Y){\PNt}{\xi}\\
&+{\PNP}AY-{\eta}(AY){\KN}+{\EKN}AY\}-g({\PY},{\phi}X){\ENK}{\KN}\Big]\\
&+(Yh)A{\phi}X+h({\Na}_YA){\phi}X-({\Na}_YA^2){\phi}X.
\end{split}
\end{equation}

The second term of \eqref{3.6} becomes
\begin{equation}\label{3.8}
\begin{split}
S(&{\Na}_Y{\phi})X=S\{{\eta}(X)AY-g(AY,X){\xi}\}\\
=&{\eta}(X)\Big[-\frac{1}{2}\Big\{(4m+7)AY-3{\eta}(AY){\xi}-3{\SN}{\EN}(AY){\KN}\\
&+{\SN}\{{\ENK}{\PNP}AY-{\eta}({\PN}AY){\PNK}-{\eta}(AY){\ENK}{\KN}\}\Big\}\\
&+hA^2Y-A^3Y\Big]\\
&-g(AY,X)\Big[-\frac{1}{2}\Big\{(4m+7){\xi}-3{\xi}-4{\SN}{\ENK}{\KN}\Big\}+hA{\xi}-A^2{\xi}\Big].
\end{split}
\end{equation}
The first term of the right side in \eqref{3.6} becomes
$$
({\Na}_Y{\phi})SX={\eta}(SX)AY-g(AY,SX){\xi},$$
and the second term of the right side in \eqref{3.6} is given by
\begin{equation}\label{3.9}
\begin{split}
{\phi}&({\nabla}_YS)X=\frac{3}{2}{\eta}(X){\phi}^2AY\\
&+\frac{3}{2}{\SN}\{{\qNt}(Y){\ENo}(X)-{\qNo}(Y){\ENt}(X)+g({\PN}AY,{\phi}X)\}{\PKN}\\
&+\frac{3}{2}{\SN}{\EN}(X)\{{\qNt}(Y){\PKNo}-{\qNo}(Y){\PKNt}+{\PPN}AY\}\\
&-\frac{1}{2}{\SN}\Big[ Y({\EN}({\xi})){\PPNP}X +{\EN}({\xi})\{-{\qNo}(Y){\phi}{\PNt}{\phi}X+{\qNt}(Y){\phi}{\PNoP}X\\
&+{\EN}({\phi}X){\phi}AY-g(AY,{\phi}X){\PKN}\}\\
&+{\EN}({\xi})\{{\eta}(X){\PPN}AY - g(AY,X){\PPNK}\}-g({\phi}AY,{\PN}X){\phi}{\PNK}\\
&+\{{\qNo}(Y){\eta}({\PNt}X)-{\qNt}(Y){\eta}({\PNo}X)\\
&-{\EN}(X){\eta}(AY)+{\EKN}g(AY,X)\}{\phi}{\PNK}\\
&-{\eta}({\PN}X)\{{\qNt}(Y){\phi}{\PNo}{\xi}-{\qNo}(Y){\PPNtK}+{\PPNP}AY\\
&-{\eta}(AY){\phi}{\KN}+{\EKN}{\phi}AY\}\\
&-g({\PY},X){\ENK}{\PKN}-{\eta}(X)Y({\ENK}){\PKN}-{\eta}(X){\ENK}{\phi}{\Na}_Y{\KN}\Big]\\
&+(Yh){\phi}AX+h{\phi}({\Na}_YA)X-{\phi}({\Na}_YA^2)X.
\end{split}
\end{equation}

Putting $X={\xi}$ into \eqref{3.5} and using that the Reeb vector field $\xi$ is principal, that is, $A{\xi}={\alpha}{\xi}$, then
we have
\begin{equation}\label{3.10}
\begin{split}
S(&{\Na}_Y{\phi}){\xi}=\Big[ -\frac{1}{2}\Big\{(4m+7)AY-3{\eta}(AY){\xi}-3{\SN}{\EN}(AY){\KN}\\
&+{\SN}\{{\ENK}{\PNP}AY-{\eta}({\PNP}AY){\PNK}-{\alpha}{\eta}(Y){\ENK}{\KN}\}\Big\}+hA^2Y-A^3Y\Big]\\
&-{\alpha}{\eta}(Y)\Big[ -\frac{1}{2}\Big\{4(m+1){\xi}-4{\SN}{\ENK}{\KN}\Big\}+({\alpha}h-{\alpha}^2){\xi}\Big] .
\end{split}
\end{equation}
Moreover, the right side of \eqref{3.6} becomes
\begin{equation}\label{3.11}
\begin{split}
({\nabla}_Y&{\phi})S{\xi}+{\phi}({\Na}_YS){\xi}\\
=& {\eta}(S{\xi})AY-g(AY,S{\xi}){\xi}+{\phi}({\Na}_YS){\xi}\\
=&\Big[\{-2(m+1)+h{\alpha}-{\alpha}^2\}+2{\SN}{\ENK}^2\Big]AY+\frac{3}{2}{\phi}^2AY\\
&-\Big[\{-2(m+1){\alpha}+h{\alpha}^2-{\alpha}^3\}{\eta}(Y)+2{\SN}{\ENK}{\EN}(AY)\Big]{\xi}\\
&+\frac{3}{2}{\SN}\{{\qNt}(Y){\ENo}({\xi})-{\qNo}(Y){\ENt}({\xi})+{\EN}({\phi}AY)\}{\PKN}\\
&+\frac{3}{2}{\SN}{\ENK}\{{\qNt}(Y){\PKNo}-{\qNo}(Y){\PKNt}+{\PPN}AY\}\\
&-\frac{1}{2}{\SN}\Big[{\ENK}\{{\PPN}AY-{\alpha}{\eta}(Y){\phi}^2{\KN}\}-g({\phi}AY,{\PKN}){\phi}^2{\KN}\\
&-Y({\ENK}){\PKN}-{\ENK}{\phi}{\Na}_Y{\KN}\Big] + h{\phi}({\Na}_YA){\xi}-{\phi}({\Na}_YA^2){\xi}.
\end{split}
\end{equation}
From this, putting $Y={\xi}$ into \eqref{3.11}, we obtain
\begin{equation}\label{3.12}
\begin{split}
0=&{\SN}\{{\qNt}({\xi}){\ENoK}-{\qNoK}{\ENtK}\}{\PKN}\\
&+{\SN}{\ENK}\{{\qNtK}{\PKNo}-{\qNoK}{\PKNt}+{\alpha}{\phi}^2{\KN}\}.
\end{split}
\end{equation}
Now in order to show that $\xi$ belongs to either the distribution $\mathcal Q$ or to the distribution ${\mathcal Q}^{\bot}$, let us assume that
${\xi}=X_1+X_2$ for some $X_1{\in}{\mathcal Q}$ and $X_2{\in}{\mathcal Q}^{\bot}$. Then it follows that

\begin{equation}\label{3.13}
\begin{split}
0=&{\SN}\Big\{{\qNt}({\xi}){\ENo}({\xi})-{\qNo}({\xi}){\ENt}({\xi})\Big\}({\PN}X_1+{\PN}X_2)\\
&+{\SN}{\ENK}\Big\{{\qNt}({\xi})({\PNo}X_1+{\PNo}X_2)\\
&-{\qNo}({\xi})({\PNt}X_1+{\PNt}X_2)-{\alpha}{\KN}+{\alpha}{\EKN}(X_1+X_2)\Big\}.
\end{split}
\end{equation}
Then by comparing $\mathcal Q$ and ${\mathcal Q}^{\bot}$ component of \eqref{3.13}, we have respectively
\begin{equation}\label{3.14}
\begin{split}
0=&{\SN}\{{\qNtK}{\ENoK}-{\qNoK}{\ENtK}\}{\PN}X_1+{\alpha}{\SN}{\ENK}^2X_1\\
&+{\SN}{\ENK}\{{\qNtK}{\PNo}X_1-{\qNoK}{\PNt}X_1\},
\end{split}
\end{equation}

\begin{equation}\label{3.15}
\begin{split}
0=&{\SN}\{{\qNtK}{\ENoK}-{\qNoK}{\ENtK}\}{\PN}X_2\\
&+{\SN}{\ENK}\{{\qNtK}{\PNo}X_2-{\qNoK}{\PNt}X_2-{\alpha}{\KN}+{\alpha}{\EKN}X_2\}.
\end{split}
\end{equation}
Taking an inner product \eqref{3.14} with $X_1$, we have
\begin{equation}\label{3.16}
{\alpha}{\SN}{\ENK}^2=0.
\end{equation}
Then ${\alpha}=0$ or ${\ENK}=0$ for ${\nu}=1,2,3$. So for a non-vanishing geodesic Reeb flow we have
${\ENK}=0$, ${\nu}=1,2,3.$ This means that ${\xi}{\in}{\mathcal Q}$, which contradicts to our assumption ${\xi}=X_1+X_2$.
Including this one, we are able to assert the following:
\vskip 6pt
\par
\begin{lemma}\quad Let $M$ be a Hopf hypersurface in $\NBt$ with commuting Ricci tensor. Then the Reeb vector
$\xi$ belongs to either the distribution ${\mathcal Q}$ or to the distribution ${\mathcal Q}^{\bot}$.
\end{lemma}
\par
\begin{proof}\quad When the geodesic Reeb flow is non-vanishing, that is ${\alpha}{\not =}0$, \eqref{3.16} gives ${\xi}{\in}{\mathcal Q}$.
When the geodesic Reeb flow is vanishing, we differentiate $A{\xi}=0$. Then by Suh (\cite{S1} and \cite{S2}) we know that
$${\SN}{\ENK}{\EN}({\phi}Y)=0.$$
From this, by replacing $Y{\in}{\mathcal Q}$ by ${\phi}Y$, it follows that
$$
{\SN}{\ENSK}{\eta}(Y)=0.$$
So if there are some $Y{\in}{\mathcal Q}$ such that ${\eta}(Y){\not =}0$, then ${\ENK}=0$ for ${\nu}=1,2,3.$
This means that ${\xi}{\in}{\mathcal Q}$. If ${\eta}(Y)=0$ for any $Y{\in}{\mathcal Q}$, then we know ${\xi}{\in}{\mathcal Q}^{\bot}$.
\end{proof}
\vskip 8pt

\section{Real hypersurfaces with geodesic Reeb flow satisfying ${\xi}{\in}{\mathcal Q}$}\label{section 4}
\setcounter{equation}{0}
\renewcommand{\theequation}{4.\arabic{equation}}
\vspace{0.13in}

Now in this section let us show that the distribution $\mathcal Q$ of a
Hopf real hypersurface $M$ in $\NBt$ with ${\xi}{\in}{\mathcal Q}$ satisfies
$g(A{\mathcal Q},{\mathcal Q}^{\bot})=0$.
\par
\vskip 6pt

The Reeb vector $\xi$ is said to be a {\it Hopf} vector if it is a
principal vector for the shape operator $A$ of $M$ in $\NBt$, that
is, the Reeb vector $\xi$ is invariant under the shape operator $A$.
\vskip 6pt
\par
In a theorem due to Berndt and Suh \cite{BS3} we know that the
Reeb vector $\xi$ of $M$ belongs to the maximal quaternionic subbundle ${\mathcal Q}$ when
$M$ is locally congruent to a real hypersurface of type $(B)$, that is,
a tube over a totally real totally geodesic ${\Bbb H}P^n$, $m=2n$, or
a horosphere in $SU_{2,m}/S(U_2U_m)$ with $JN{\bot}{\frak J}N$ whose center at infinity is singular, and a real hypersurface of type $(D)$ in $\NBt$.  Naturally we are
able to consider a converse problem. From such a view point Suh \cite{S6} has proved the following
for real hypersurfaces in $\NBt$ with the Reeb vector field $\xi{\in}\mathcal Q$.
\vskip 6pt
\par
\begin{thm C}\label{Btype}
Let $M$ be a real hypersurface in noncompact
complex two-plane Grassmannian $\NBt$ with the
Reeb vector field belonging to the maximal quaternionic subbundle $\Cal Q$. Then one
of the following statements holds,
\par
$(B)$\quad $M$ is an open part of a tube around a totally geodesic
${\Bbb H}H^n$ in \hfill
\newline
$SU_{2,2n}/S(U_2U_{2n})$, $m = 2n$,
\par
$(C_2)$\quad $M$ is an open part of a horosphere in
$SU_{2,m}/S(U_2U_m)$ whose center at infinity is singular and of
type $JN \perp {\frak J}N$,
\par
or the following exceptional case holds:
\par
$(D)$\quad The normal bundle $\nu M$ of $M$ consists of singular
tangent vectors of type $JX \perp {\frak J}X$. Moreover, $M$ has at
least four distinct principal curvatures, three of which are given
by
\begin{equation*}
\alpha = \sqrt{2}\ ,\ \gamma = 0\ ,\ \lambda = \frac{1}{\sqrt{2}}
\end{equation*}
with corresponding principal curvature spaces
\begin{equation*}
T_\alpha = TM \ominus ({\mathcal C} \cap {\mathcal Q})\ ,\ T_\gamma = J(TM
\ominus {\mathcal Q})\ ,\ T_\lambda \subset  {\mathcal C} \cap {\mathcal Q} \cap
J{\mathcal Q}.
\end{equation*}
If $\mu$ is another (possibly nonconstant) principal curvature
function, then we have $T_\mu \subset {\mathcal C} \cap {\mathcal Q} \cap
J{\mathcal Q}$, $JT_\mu \subset T_\lambda$ and ${\frak J}T_\mu \subset
T_\lambda $.
\end{thm C}
\vskip 6pt
\par
In the proof of Theorem C we have used the equation of Codazzi in section 2 for Hopf real hypersurfaces $M$ in $\NBt$, and proved that the quaternionic maximal subbundle $\mathcal Q$ is invariant
under the shape operator, that is, $g(A{\mathcal Q},{\mathcal Q}^{\bot})=0$, if the Reeb vector field $\xi$ belongs to the subbundle $\mathcal Q$ of $M$. So by using a theorem due to Berndt and Suh \cite{BS3}
we can assert Theorem B in the introduction. Then among the classification of Theorem A the Reeb vector field $\xi$ of real hypersurfaces in Theorem C belongs to the maximal quaternionic subbundle $\mathcal Q$.
\par
\vskip 6pt

\vskip 6pt
\par
Now let us check whether the Ricci tensor $S$ of hypersurfaces mentioned in Theorem C satisfies the commuting condition or not. In order to do this, we should find all of the principal curvatures
corresponding to the hypersurfaces in Theorem B.  For cases of type $(B)$, one of type $(C)$
which will be said to be of type $(C_2)$, and of type $(D)$ in
Theorem B let us introduce a proposition due to Berndt and Suh [5] as follows:
\vskip 6pt
\begin{proposition}\label{Proposition B}
\quad Let $M$ be a connected hypersurface
in $SU_{2,m}/S(U_2U_m)$, $m \geq 2$. Assume that the maximal complex
subbundle $\mathcal {C}$ of $TM$ and the maximal quaternionic subbundle
${\mathcal Q}$ of $TM$ are both invariant under the shape operator of
$M$. If $JN \perp {\frak J}N$, then one of the following statements
holds:
\par
$(B)$\quad $M$ has five (four for $r = \sqrt{2}{\text Arctanh}(1/\sqrt{3})$ in which case $\alpha = \lambda_2$)
distinct constant principal curvatures
\begin{equation*}
\begin{split}
\alpha = &\sqrt{2}\tanh(\sqrt{2}r)\ ,\ \beta =
\sqrt{2}\coth(\sqrt{2}r)\ ,\ \gamma = 0\ , \\
\lambda_1 =& \frac{1}{\sqrt{2}}\tanh(\frac{1}{\sqrt{2}}r)\ ,\ \lambda_2 =
\frac{1}{\sqrt{2}}\coth(\frac{1}{\sqrt{2}}r),
\end{split}
\end{equation*}
and the corresponding principal curvature spaces are
\begin{equation*}
T_\alpha = TM \ominus \mathcal {C}\ ,\ T_\beta = TM \ominus \mathcal {Q}\ ,\ T_\gamma =
J(TM \ominus \mathcal {Q}) = JT_\beta.
\end{equation*}
The principal curvature spaces $T_{\lambda_1}$ and $T_{\lambda_2}$
are invariant under ${\frak J}$ and are mapped onto each other
by $J$. In particular, the quaternionic dimension of
$SU_{2,m}/S(U_2U_m)$ must be even.
\par
$(C_2)$\quad $M$ has exactly three distinct constant principal curvatures
\begin{equation*}
\alpha = \beta = \sqrt{2}\ ,\ \gamma = 0\ ,\ \lambda =
\frac{1}{\sqrt{2}}
\end{equation*}
with corresponding principal curvature spaces
\begin{equation*}
T_\alpha = TM \ominus (\mathcal {C} \cap \mathcal {Q})\ ,\ T_\gamma =
J(TM \ominus \mathcal {Q})\ ,\ T_\lambda = \mathcal {C} \cap \mathcal {Q} \cap J\mathcal {Q}.
\end{equation*}
\par
$(D)$\quad $M$ has at least four distinct principal curvatures,
three of which are given by
\begin{equation*}
\alpha = \beta = \sqrt{2}\ ,\ \gamma = 0\ ,\ \lambda =
\frac{1}{\sqrt{2}}
\end{equation*}
with corresponding principal curvature spaces
\begin{equation*}
T_\alpha = TM \ominus (\mathcal {C} \cap \mathcal {Q})\ ,\ T_\gamma =
J(TM \ominus {\mathcal Q})\ ,\ T_\lambda \subset \mathcal {C} \cap \mathcal {Q} \cap J\mathcal {Q}.
\end{equation*}
If $\mu$ is another (possibly nonconstant) principal
curvature function  , then $JT_{\mu} \subset T_\lambda$
and ${\frak J}T_{\mu} \subset T_\lambda$.
\end{proposition}
\vskip 6pt
On the other hand, we calculate the following
\begin{equation*}
\begin{split}
S{\phi}X=&-\frac{1}{2}\Big[(4m+7){\phi}X-3{\SN}{\EN}({\phi}X){\KN}+{\SN}\{{\ENK}{\PNP}^2X-{\eta}({\PNP}X){\PNK}\\
&-{\eta}({\phi}X){\ENK}){\KN}\}\Big]+hA{\phi}X-A^2{\phi}X
\end{split}
\end{equation*}
and
\begin{equation*}
\begin{split}
{\phi}SX=&-\frac{1}{2}\Big[(4m+7){\phi}X-3{\SN}{\EN}(X){\PKN}+{\SN}\{{\ENK}{\PPNP}X-{\eta}({\PN}X){\PPNK}\\
&-{\eta}(X){\ENK}{\PKN}\}\Big]+h{\phi}AX-{\phi}A^2X.
\end{split}
\end{equation*}
Then the Ricci commuting $S{\phi}={\phi}S$ implies that
\begin{equation}\label{4.1}
hA{\phi}X-A^2{\phi}X=h{\phi}AX-{\phi}A^2X+2{\SN}{\EN}(X){\PKN}-2{\SN}{\EN}({\phi}X){\KN}.
\end{equation}

First let us check that real hypersurfaces of type $(B)$ in Theorem B satisfy the formula \eqref{4.1} or not.
\par
Putting $X={\xi}{\in}{\mathcal Q}$ we know the formula (4.1) holds. For $X={\xi}_i{\in}{\mathcal Q}^{\bot}$
the left side of \eqref{4.1} becomes
$$L=hA{\phi}{\xi}_i-A^2{\phi}{\xi}_i=0$$
and the right side is given by
\begin{equation*}
\begin{split}
R=&h{\phi}A{\xi}_i-{\phi}A^2{\xi}_i+2{\EN}({\xi}_i){\PKN}-2{\SN}{\EN}({\phi}{\xi}_i){\KN}\\
=&h{\beta}{\phi}{\xi}_i-{\beta}^2{\phi}{\xi}_i+2{\phi}{\xi}_i=(h{\beta}-{\beta}^2+2){\phi}{\xi}_i .
\end{split}
\end{equation*}
Then ${\beta}(h-{\beta})=-2$.
\par
\vskip 6pt
On the other hand, the trace $h$ is given by
\begin{equation*}
\begin{split}
h=&{\alpha}+3{\beta}+3{\gamma}+2(m-2)({\lambda_1}+{\lambda_2})\\
=&{\sqrt 2}\tanh ({\sqrt 2}r)+(2m-1){\sqrt 2}\coth ({\sqrt 2}r),
\end{split}
\end{equation*}
where we have used $2\coth {\sqrt 2}r = \tanh (\frac{1}{\sqrt 2}r)+ \coth (\frac{1}{\sqrt 2}r)$.
Then it follows that
\begin{equation*}
\begin{split}
{\beta}(h-{\beta})=&({\sqrt 2}\tanh ({\sqrt 2}r)+(2m-2){\sqrt 2}\coth ({\sqrt 2}r)){\sqrt 2}\coth ({\sqrt 2}r)\\
=&2 + 4(m-1)\coth^2 ({\sqrt 2}r)\\
=&-2,
\end{split}
\end{equation*}
which gives $4(m-1)\coth^2 ({\sqrt 2}r)=-4$. This gives a contradiction.
\par
\vskip 6pt
As a second,\ let us check whether or not a horosphere at infinity could satisfy the Ricci commuting.  By Proposition 4.1, we know that
the principal curvatures of the case $(C_2)$ are given by
$${\alpha}={\beta}={\sqrt 2}, {\gamma}=0, {\lambda}=\frac{1}{\sqrt 2}.$$
Then the trace of the shape operator $h$ becomes
\begin{equation*}
\begin{split}
h=&{\alpha}+3{\beta}+3{\gamma}+4(m-2){\lambda}\\
=&4{\sqrt 2}+2(m-2){\sqrt 2}
\end{split}
\end{equation*}
Then by putting $X={\xi}_i{\in}{\mathcal Q}^{\bot}$ in (4.1), we know that $(h{\beta}-{\beta}^2+2){\phi}{\xi}_i=0$. This gives that
$-2={\beta}(h-{\beta})=2(2m-1)$, which gives a contradiction.
\par
\vskip 6pt
As a third,\ we want to check that real hypersurfaces of type $(D)$ in Theorem C satisfy the Ricci commuting or not. Then the principal curvature of type $(D)$ becomes
$${\alpha}={\beta}={\sqrt 2}, {\gamma}=0, {\lambda}=\frac{1}{\sqrt 2}\quad \text{and}\ quad \mu $$
where $JT_{\mu}{\subset}T_{\lambda}$ and ${\frak J}T_{\mu}{\subset}T_{\lambda}$. Then the trace $h$ becomes that
\begin{equation*}
\begin{split}
h=&{\alpha}+3{\beta}+3{\gamma}+(2m-4)\frac{1}{\sqrt 2}+(2m-4){\mu}\\
=&(m-2)({\sqrt 2}+2{\mu})+4{\sqrt 2}.
\end{split}
\end{equation*}
From this, together with ${\beta}(h-{\beta})=-2$, we have
$${\sqrt 2}\{(m-2)({\sqrt 2}+2{\mu})+3{\sqrt 2}\}=-2.$$
This gives $2{\mu}(m-2)=-(m+2){\sqrt 2}$. Then the trace $h$ becomes that
\begin{equation*}
\begin{split}
h=&(m-2)({\sqrt 2}+2{\mu})+4{\sqrt 2}\\
=&(m-2){\sqrt 2}-(m+2){\sqrt 2}+4{\sqrt 2}\\
=&0.
\end{split}
\end{equation*}
On the other hand, from the Ricci commuting (4.1) we know that
$$hA{\phi}X-A^2{\phi}X=h{\phi}AX-{\phi}A^2X$$
for any $X{\in}{\mathcal Q}$. So it follows that $h={\lambda}+{\mu}=-\frac{2{\sqrt 2}}{m-2}$. But this is a contradiction.
So also the case of type $(D)$ can not occur.
\par
\vskip 6pt
Summing up all cases $(B)$, $(C_2)$ and $(D)$ mentioned above, we conclude that there do not exist any Hopf real hypersurfaces in complex hyperbolic two-plane Grassmannians
$\NBt$ with commuting Ricci tensor when the Reeb vector field $\xi$ belongs to the distribution $\mathcal Q$.
\par
\vskip 6pt
\section{Real hypersurfaces with geodesic Reeb flow satisfying ${\xi}{\in}{\mathcal Q}^{\bot}$}\label{section 5}
\setcounter{equation}{0}
\renewcommand{\theequation}{5.\arabic{equation}}
\vskip 6pt
\par
Now let us consider a Hopf real hypersurface $M$ in $\NBt$ with commuting Ricci tensor satisfying ${\xi}{\in}{\mathcal Q}^{\bot}$, where ${\mathcal Q}$ denotes a
quaternionic maximal subbundle in $T_xM$, $x{\in}M$ such that $T_xM={\mathcal Q}{\oplus}{\mathcal Q}^{\bot}$.
Since we assume $\xi{\in}{\mathcal Q}^{\bot}= \text{Span}\
\{{\Ko},{\Kt},{\Ks}\}$, there exists a Hermitian structure
$J_1{\in}{\frak J}$ such that $JN=J_1N$, that is, ${\xi}={\xi}_1$.
\par
\vskip 6pt
Moreover, the right side of \eqref{3.6} can be written as follows:
\begin{equation}\label{5.1}
\begin{split}
({\nabla}_Y&{\phi})S{\xi}+{\phi}({\Na}_YS){\xi}\\
=& {\eta}(S{\xi})AY-g(AY,S{\xi}){\xi}+{\phi}({\Na}_YS){\xi}\\
=&\Big[ \{-2(m+1)+h{\alpha}-{\alpha}^2\}+2{\SN}{\ENK}^2\Big]AY+\frac{3}{2}{\phi}^2AY\\
&-\Big[\{-2(m+1){\alpha}+h{\alpha}^2-{\alpha}^3\}{\eta}(Y)+2{\SN}{\ENK}{\EN}(AY)\Big]{\xi}\\
&+\frac{3}{2}{\SN}\{{\qNt}(Y){\ENo}({\xi})-{\qNo}(Y){\ENt}({\xi})+{\EN}({\phi}AY)\}{\PKN}\\
&+\frac{3}{2}{\SN}{\ENK}\{{\qNt}(Y){\PKNo}-{\qNo}(Y){\PKNt}+{\PPN}AY\}\\
&-\frac{1}{2}{\SN}\Big[{\ENK}\{{\PPN}AY-{\alpha}{\eta}(Y){\phi}^2{\KN}\}-g({\phi}AY,{\PKN}){\phi}^2{\KN}\\
&-Y({\ENK}){\PKN}-{\ENK}{\phi}{\Na}_Y{\KN}\Big] + h{\phi}({\Na}_YA){\xi}-{\phi}({\Na}_YA^2){\xi}.
\end{split}
\end{equation}
Then $({\Na}_YS){\phi}{\xi}+S({\Na}_Y{\phi}){\xi}=({\Na}_Y{\phi})S{\xi}+{\phi}({\Na}_YS){\xi}$ and ${\xi}{\in}{\mathcal Q}^{\bot}$,
that is, ${\xi}={\xi}_1$, becomes
\begin{equation}\label{5.2}
\begin{split}
-\frac{1}{2}&\Big\{(4m+7)AY-6{\alpha}{\eta}(Y){\xi}-3{\Et}(AY){\Kt}-3{\Es}(AY){\Ks}\\
&+{\PoP}AY-{\eta}({\Pt}AY){\PtK}-{\eta}({\Ps}AY){\PsK}-{\alpha}{\eta}(Y){\xi}\Big\}\\
&+hA^2Y-A^3Y-{\alpha}{\eta}(Y)\{-2m{\xi}+({\alpha}h-{\alpha}^2){\xi}\}\\
=&\{-2m+h{\alpha}-{\alpha}^2\}AY-\{-2m+h{\alpha}-{\alpha}^2\}AY\\
&-\{-2m+h{\alpha}-{\alpha}^2\}{\alpha}{\eta}(Y){\xi}\\
\end{split}
\end{equation}
\begin{equation*}
\begin{split}
&+\frac{3}{2}\Big[{\phi}^2AY+(q_1(Y){\EsK}-q_3(Y){\eta}({\xi})\\
&+{\Et}({\phi}AY){\PKt}){\PKt}+(q_2(Y){\EK}-q_1(Y){\EtK}+{\Es}({\phi}AY){\PKs}){\PKs}\\
&+(q_3(Y){\PKt}-q_2(Y){\PKs}+{\PPo}AY)\Big]\\
&-\frac{1}{2}\Big[{\PPo}AY-{\alpha}{\eta}(Y){\PPoK}-g({\phi}AY,{\PtK}){\PPtK}-g({\phi}AY,{\PsK}){\PPsK}\\
&-{\phi}^2AY\Big]+h{\alpha}{\phi}^2AY-h{\phi}A{\phi}AY-{\alpha}^2{\phi}^2AY+{\phi}A^2{\phi}AY.
\end{split}
\end{equation*}
Then \eqref{5.2} can be rearranged as follows:
\begin{equation}\label{5.3}
\begin{split}
-5AY&+5{\alpha}{\eta}(Y){\xi}+6{\Et}(AY){\Kt}+6{\Es}(AY){\Ks}+2hA^2Y-2A^3Y\\
=&3{\PPo}AY+2(h{\alpha}^2-{\alpha}^3){\eta}(Y){\xi}-2h{\phi}A{\phi}AY+2{\phi}A^2{\phi}AY.
\end{split}
\end{equation}

\vskip 6pt
Now let us use a Proposition due to Suh (see \cite{S5}) as follows:

\begin{proposition}\label{Proposition 5.1}
\quad If $M$ is a connected orientable real hypersurface in complex
hyperbolic two-plane Grassmannian $\NBt$ with geodesic Reeb flow,
then
\begin{equation*}
\begin{split}
2g&(A\phi AX,Y) - \alpha g((A\phi + \phi A)X,Y) +  g(\phi X,Y)\\
=& {\SN} \Big\{{\EN}(X){\EN}({\phi}Y) - {\EN}(Y){\EN}({\phi}X) - g({\PN}X,Y){\ENK}\\
&-2{\eta}(X){\EN}({\phi}Y){\ENK} + 2 {\eta}(Y){\EN}({\phi}X){\ENK}\Big\}\ .
\end{split}
\end{equation*}
\end{proposition}

Then by Proposition 5.1, we know that for ${\xi}={\xi}_1$
\begin{equation}\label{5.4}
2A{\phi}AY={\alpha}(A{\phi}+{\phi}A)Y-{\phi}Y-{\phi}_1Y+2\{{\Et}(Y){\Ks}-{\Es}(Y){\Kt}\}.
\end{equation}
So it follows that $2A{\phi}AY={\alpha}(A{\phi}+{\phi}A)Y-{\phi}Y-{\phi}_1Y$ for any $Y{\in}{\mathcal Q}$, where ${\mathcal Q}$ denotes a maximal quaternionic subbundle of $T_xM$, $x{\in}M$. Then this gives the following
\begin{equation*}
\begin{split}
2{\phi}A^2{\phi}AY=&2{\phi}A(A{\phi}AY)={\phi}A\{{\alpha}(A{\phi}+{\phi}A)Y-{\phi}Y-{\Po}Y\}\\
=&{\alpha}{\phi}A(A{\phi}+{\phi}A)Y-{\phi}A{\phi}Y-{\phi}A{\Po}Y.
\end{split}
\end{equation*}
From this, together with \eqref{5.4}, it follows that
\begin{equation}\label{5.5}
\begin{split}
-5AY&+5{\alpha}{\eta}(Y){\xi}+6{\Et}(AY){\Kt}+6{\Es}(AY){\Ks}+2hA^2Y-2A^3Y\\
=&3{\PPo}AY+2(h{\alpha}-{\alpha}^2){\alpha}{\eta}(Y){\xi}-h{\phi}\{{\alpha}(A{\phi}+{\phi}A)Y-{\phi}Y-{\Po}Y\}\\
&+{\alpha}{\phi}A(A{\phi}+{\phi}A)Y-{\phi}A{\phi}Y-{\phi}A{\Po}Y
\end{split}
\end{equation}

Then \eqref{5.5} can be written as follows:
\begin{equation}\label{5.6}
\begin{split}
3{\PPo}&AY-2hA^2Y+2A^3Y+5AY=h{\alpha}{\phi}(A{\phi}+{\phi}A)Y+hY-h{\PPo}Y\\
&+6{\Et}(AY){\Kt}+6{\Es}(AY){\Ks}-{\alpha}{\phi}A(A{\phi}+{\phi}A)Y\\
&+{\phi}A{\phi}Y+{\phi}A{\Po}Y.
\end{split}
\end{equation}

On the other hand, the commuting Ricci tensor ${\phi}S=S{\phi}$ gives for any $Y{\in}{\mathcal Q}$ and ${\xi}={\xi}_1{\in}{\mathcal Q}^{\bot}$ 
\begin{equation*}
\begin{split}
-\frac{1}{2}&\Big[(4m+7){\phi}Y+{\PoP}^2Y\Big]+hA{\phi}Y-A^2{\phi}Y\\
=&-\frac{1}{2}\Big[(4m+7){\phi}Y+{\PPo}{\phi}Y\Big]=h{\phi}AY-{\phi}A^2Y
\end{split}
\end{equation*}
Then for any $Y$ belonging to the maximal quaternionic subbundle $\mathcal Q$, we know
$$hA{\phi}Y-A^2{\phi}Y=h{\phi}AY-{\phi}A^2Y.$$
From this, by replacing $Y$ by ${\phi}Y$ for $Y{\in}{\mathcal Q}$ and applying $A$ to the obtained equation, we have
\begin{equation}\label{5.7}
hA^2Y-A^3Y=-hA{\phi}A{\phi}Y+A{\phi}A^2{\phi}Y.
\end{equation}
\par
\vskip 6pt
On the other hand, by using \eqref{5.4} for any ${\phi}Y{\in}{\mathcal Q}$, we know
$$2A{\phi}A{\phi}Y={\alpha}(A{\phi}+{\phi}A){\phi}Y+Y-{\PoP}Y$$
and
\begin{equation*}
\begin{split}
2A{\phi}A^2&{\phi}Y={\alpha}(A{\phi}+{\phi}A)A{\phi}Y-{\phi}A{\phi}Y-{\Po}A{\phi}Y\\
&+2\{{\Et}(A{\phi}Y){\Ks}-{\Es}(A{\phi}Y){\Ks}\}.
\end{split}
\end{equation*}
Then substituting these formulas into \eqref{5.7}, we have
\begin{equation*}
\begin{split}
-2h&A^2Y+2A^3Y=2hA{\phi}A{\phi}Y-2A{\phi}A^2{\phi}Y\\
=&h{\alpha}(A{\phi}+{\phi}A){\phi}Y+hY-h{\PoP}Y\\
&-{\alpha}(A{\phi}+{\phi}A)A{\phi}Y+{\phi}A{\phi}Y\\
&-2\{{\Et}(A{\phi}Y){\Ks}-{\Es}(A{\phi}Y){\Kt}\}.
\end{split}
\end{equation*}
From this, together with \eqref{5.6}, it follows that
\begin{equation}\label{5.8}
\begin{split}
3{\PoP}&AY+5AY+ h{\alpha}(A{\phi}+{\phi}A){\phi}Y+hY-h{\PoP}Y\\
&-{\alpha}(A{\phi}+{\phi}A)A{\phi}Y+{\phi}A{\phi}Y+{\phi}_1A{\phi}Y\\
&-2\{{\Et}(A{\phi}Y){\Ks}-{\Es}(A{\phi}Y){\Kt}\}\\
=&{\alpha}h(A{\phi}+{\phi}A)+hY-h{\PPo}Y+6{\Et}(AY){\Kt}+6{\Es}(AY){\Ks}\\
&-{\alpha}{\phi}A(A{\phi}+{\phi}A)Y+{\phi}A{\phi}Y+{\phi}A{\phi}_1Y.
\end{split}
\end{equation}
Then it can be rearranged as follows:
\begin{equation}\label{5.9}
\begin{split}
6{\PPo}AY&+10AY-2h{\alpha}AY-2{\alpha}A{\phi}A{\phi}Y+2{\Po}A{\phi}Y\\
&-4\{{\Et}(A{\phi}Y){\Ks}-{\Es}(A{\phi}Y){\Kt}\}\\
=&-2{\alpha}hAY+12{\Et}(AY){\Kt}+12{\Es}(AY){\Ks}\\
&-2{\alpha}{\phi}A{\phi}AY+2{\phi}A{\Po}Y.
\end{split}
\end{equation}

On the other hand, from (5.4) we calculate the following formulas
$$
2{\alpha}A{\phi}A{\phi}Y={\alpha}^2(A{\phi}+{\phi}A){\phi}Y+{\alpha}Y-{\alpha}{\phi}_1{\phi}Y$$
and
$$
2{\alpha}{\phi}A{\phi}AY={\alpha}^2{\phi}(A{\phi}+{\phi}A)Y+{\alpha}Y-{\alpha}{\phi}{\phi}_1Y.$$
Substituting these formulas into \eqref{5.9}, we have
\begin{equation}\label{5.10}
\begin{split}
&6{\PPo}AY+10AY-2h{\alpha}AY-{\alpha}^2(A{\phi}+{\phi}A){\phi}Y-{\alpha}Y+{\alpha}{\PoP}Y+2{\Po}A{\phi}Y\\
&-4\{{\Et}(A{\phi}Y){\Ks}-{\Es}(A{\phi}Y){\Kt}\}\\
=&-2{\alpha}hAY+12{\Et}(AY){\Kt}+12{\Es}(AY){\Ks}-{\alpha}^2{\phi}(A{\phi}+{\phi}A)Y\\
&-{\alpha}Y+{\alpha}{\PPo}Y+2{\phi}A{\Po}Y.
\end{split}
\end{equation}
Thus it can be rearranged for any $Y{\in}{\mathcal Q}$:
\begin{equation}\label{5.11}
\begin{split}
&6{\PPo}AY+10AY+2{\Po}A{\phi}Y-4\{{\Et}(A{\phi}Y){\Ks}-{\Es}(A{\phi}Y){\Kt}\}\\
=&12{\Et}(AY){\Kt}+12{\Es}(AY){\Ks}+{\alpha}^2AY+2{\phi}A{\Po}Y.
\end{split}
\end{equation}
From this, if we take an inner product with $\Kt$, it follows that
$$6{\Et}(AY)+10{\Et}(AY)-2{\Es}(A{\phi}Y)+4{\Es}(A{\phi}Y)=12{\Et}(AY)+2{\Es}(A{\Po}Y).$$
Then it can be arranged by
\begin{equation}\label{5.12}
{\Es}(A{\phi}Y)=-2{\Et}(AY)+{\Es}(A{\Po}Y).
\end{equation}
Similarly, let us take an inner product \eqref{5.11} by $\Ks$. Then we have
$$6{\Es}(AY)+10{\Es}(AY)+2{\Et}(A{\phi}Y)-4{\Et}(A{\phi}Y)=12{\Es}(AY)-2{\Et}(A{\Po}Y).$$
Then it can be arranged by
\begin{equation}\label{5.13}
{\Et}(A{\phi}Y)=2{\Es}(AY)+{\Et}(A{\Po}Y)
\end{equation}
for any $Y{\in}{\mathcal Q}$.  Note that on a maximal quaternionic subbundle $\mathcal Q$ we know that
$({\PPo})^2=I$, that is, ${\PPo}{\PPo}X=X$ for any $X{\in}T_xM$, $x{\in}M$ and $\text{tr}({\PPo})=0$. By virtue of these facts, we can decompose a maximal
quaternionic subbundle $\mathcal Q$ in such two eigenspaces as follows:
$$E_{+1}=\{X{\in}{\mathcal Q}{\vert}{\PPo}X=X\}$$
and
$$E_{-1}=\{X{\in}{\mathcal Q}{\vert}{\PPo}X=-X\},$$
where ${\mathcal Q}=E_{+1}{\oplus}E_{-1}$. Then we know that
$$
X{\in}E_{+1}\ \text{iff}\ {\PPo}X=X\ \text{iff}\ {\Po}X=-{\phi}X$$
and
$$
X{\in}E_{-1}\ \text{iff}\ {\PPo}X=-X\ \text{iff}\ {\phi}X={\Po}X.$$
First let us consider on the subbundle $E_{-1}=\{Y{\in}{\mathcal Q}{\vert}{\phi}Y={\Po}Y\}$. Then \eqref{5.12} gives ${\Et}(AY)=0$ for any $Y{\in}E_{-1}$. Moreover, \eqref{5.13} implies
${\Es}(AY)=0$ for any $Y{\in}E_{-1}$.
\par
\vskip 6pt
Now let us consider a subbundle $E_{+1}=\{X{\in}{\mathcal Q}{\vert}{\phi}X=-{\Po}X\}$ in the maximal subbundle $\mathcal Q$.  Then \eqref{5.12} and \eqref{5.13} respectively gives the
following
\begin{equation}\label{5.14}
 {\Es}(A{\phi}Y)={\Et}(AY)\ \text{and}\ {\Et}(A{\phi}Y)=-{\Es}(AY).
\end{equation}
From this, together with \eqref{5.11}, we have
\begin{equation}\label{5.15}
3{\PPo}AY+5AY+{\Po}A{\phi}Y-4\{{\Et}(AY){\Kt}+{\Es}(AY){\Ks}\}={\phi}A{\Po}Y.
\end{equation}
Now let us consider eigenvectors $Y,{\phi}Y{\in}E_{+1}$ such that ${\phi}Y=-{\Po}Y$. Then we may put
$$
AY={\lambda}Y + {\SN}{\EN}(AY){\KN}
$$
$$
A{\phi}Y={\bar{\lambda}}{\phi}Y + {\SN}{\EN}(A{\phi}Y){\KN}
$$
and
$$
{\phi}AY={\lambda}{\phi}Y + {\SN}{\EN}(AY){\phi}{\KN}.
$$
From these formulas it follows that for any $Y{\in}E_{+1}$ satisfying ${\phi}Y=-{\Po}Y$
$$
{\PPo}AY={\phi}\{{\lambda}{\Po}Y+{\SN}{\EN}(AY){\Po}{\KN}\}={\lambda}Y+{\SN}{\EN}(AY){\PPo}{\KN},$$
$$
{\PoP}A{\phi}Y={\Po}\{{\bar\lambda}{\phi}Y+{\SN}{\EN}(A{\phi}Y){\KN}\}={\bar\lambda}Y+{\SN}{\EN}(A{\phi}Y){\Po}{\KN},$$
and
$$
{\phi}A{\phi}Y=-\{{\bar\lambda}Y - {\SN}{\EN}(A{\phi}Y){\PKN}\}=-{\phi}A{\Po}Y .$$
From these formulas, \eqref{5.15} gives the following
\begin{equation}\label{5.16}
\begin{split}
3\{&{\lambda}Y+{\SN}{\EN}(AY){\PPo}{\KN}\}+5\{{\lambda}Y+{\SN}{\EN}(AY){\KN}\}\\
&+{\bar{\lambda}}Y + {\SN}{\EN}(A{\phi}Y){\Po}{\KN}\\
=&4\{{\Et}(AY){\Kt}+{\Es}(AY){\Ks}\} + \{{\bar\lambda}Y - {\SN}{\EN}(A{\phi}Y){\PKN}\}
\end{split}
\end{equation}
which gives ${\lambda}=0$.  Then this implies the following
\begin{equation}\label{5.17}
AY={\SN}{\EN}(AY){\KN}=g(A{\Kt},Y){\Kt}+g(A{\Ks},Y){\Ks}.
\end{equation}
Then for any ${\phi}Y{\in}E_{+1}$ we know that
\begin{equation*}
\begin{split}
A{\phi}Y=&g(A{\Kt},{\phi}Y){\Kt}+g(A{\Ks},{\phi}Y){\Ks}\\
=&-{\Es}(AY){\Kt}+{\Et}(AY){\Ks}\\
{\phi}A{\phi}Y=&g(A{\Kt},{\phi}Y){\phi}{\Kt}+g(A{\Ks},{\phi}Y){\phi}{\Ks}\\
=&-g(A{\Kt},{\phi}Y){\Ks}+g(A{\Ks},{\phi}Y){\Kt}\\
{\phi}_1A{\phi}Y=&g(A{\Kt},{\phi}Y){\Ks}-g(A{\Ks},{\phi}Y){\Kt}.
\end{split}
\end{equation*}

Substituting these fomulas into \eqref{5.15} and using ${\phi}Y=-{\Po}Y$ for $Y{\in}E_{+1}$, we have
\begin{equation}\label{5.18}
\begin{split}
3{\PPo}AY=&-5AY-{\Po}A{\phi}Y+4\{{\Et}(AY){\Kt}+{\Es}(AY){\Ks}\}-{\phi}A{\phi}Y\\
=&-5\{{\Et}(AY){\Kt}+{\Es}(AY){\Ks}\}+4\{{\Et}(AY){\Kt}+{\Es}(AY){\Ks}\}\\
=&-\{{\Et}(AY){\Kt}+{\Es}(AY){\Ks}\}\\
=&3{\PoP}AY.
\end{split}
\end{equation}
Then by applying $\Po$ to the second equality of the above equation, we know that
$$
3{\phi}AY ={\Et}(AY){\Ks}-{\Es}(AY){\Kt}
$$
for any $Y{\in}E_{+1}$, which gives that
$$-3AY={\Et}(AY){\PKs}-{\Es}(AY){\PKt}={\Et}(AY){\Kt}+{\Es}(AY){\Ks}.$$
From this, by applying ${\Kt}$ and $\Ks$, we get the following respectively
$${\Et}(AY)=0\quad \text{and}\quad {\Es}(AY)=0$$
for any $Y{\in}E_{+1}$. Summing up the case in the subbundle $E_{-1}$ and the fact that $M$ is Hopf we proved that $g(A{\mathcal Q}, {\mathcal Q}^{\bot})=0$, that is the maximal quaternionic subbundle
$\mathcal Q$ is invariant under the shape operator.  Then from such a view point, by Theorem A in the introduction we conclude that $M$ is locally congruent to one of real hypersurfaces of type $(A)$, $(B)$, $(C_1)$, $(C_2)$ and $(D)$.
Among them we know that the Reeb vector field $\xi$ of hypersurfaces of type $(A)$ and $(C_1)$ belongs to the quaternionic maximal subbundle $\mathcal Q$, that is
$JN\in{\frak J}N$.  This gives a complete proof of our main theorem.
\qed
\par
\vskip 6pt
On the other hand, the Reeb flow of these type hypersurfaces mentioned in our theorem is isometric as in Theorem B. That is, the shape operator commutes with the structure tensor, that is $A{\phi}={\phi}A$. So by using this property, conversely, let us check whether or not these hypersurfaces satisfy the Ricci commuting.
In order to do this, let us introduce the following Proposition due to Berndt and Suh \cite{BS3} :
\par
\vskip 6pt

\begin{pro B}\label{Proposition B} \quad Let $M$ be a connected hypersurface in $SU_{2,m}/S(U_2U_m)$, $m \geq
2$. Assume that the maximal complex subbundle ${\mathcal C}$ of $TM$
and the maximal quaternionic subbundle ${\mathcal Q}$ of $TM$ are
both invariant under the shape operator of $M$. If $JN \in
{\frak J}N$, then one the following statements hold:
\par
$(A)$\quad $M$ has exactly four
distinct constant principal curvatures
\begin{equation*}
\alpha = 2\coth(2r)\ ,\ \beta = \coth(r)\ , \lambda_1 = \tanh(r)\ ,\
\lambda_2 = 0,
\end{equation*}
and the corresponding principal curvature spaces are
\begin{equation*}
T_\alpha = TM \ominus {\mathcal C}\ ,\ T_\beta = {\mathcal C} \ominus {\mathcal Q}
\ ,\ T_{\lambda_1} = E_{-1}\ ,\ T_{\lambda_2} = E_{+1}.
\end{equation*}
The principal curvature spaces $T_{\lambda_1}$ and $T_{\lambda_2}$
are complex (with respect to $J$) and totally complex (with respect
to ${\frak J}$).
\par
$(C_1)$\quad $M$ has exactly three distinct constant principal curvatures
\begin{equation*}
\alpha = 2\ ,\ \beta = 1\ ,\ \lambda = 0
\end{equation*}
with corresponding principal curvature spaces
\begin{equation*}
T_\alpha = TM \ominus {\mathcal C}\ ,\ T_\beta = ({\mathcal C} \ominus {\mathcal Q})
\oplus E_{-1}\ ,\ T_\lambda = E_{+1}.
\end{equation*}
\end{pro B}
\vskip 6pt
Now by using above proposition and the isometric Reeb flow, let us check real hypersurfaces of type $(A)$ and $(C_2)$ satisfy $S{\phi}={\phi}S$.
The from $A{\phi}={\phi}A$ the Ricci commuting gives the following
\begin{equation}\label{5.19}
\begin{split}
-3{\SN}{\EN}&({\phi}Y){\KN}+{\SN}\{{\ENK}{\PN}{\phi}^2Y-{\eta}({\PN}{\phi}Y){\PNK}\}\\
=&-3{\SN}{\EN}(Y){\PKN}+{\SN}\{{\ENK}{\PPNP}Y-{\eta}({\PN}Y){\PPNK}\}.
\end{split}
\end{equation}
\par
\vskip 6pt
Now let us check it for real hypersurfaces of type $(A)$ and of type $(C_1)$ in Proposition B as follows:
\par
\vskip 6pt
Case 1) Put $Y={\xi}={\Ko}$. Then the both sides equal to each other.
\par
Case 2) Put $Y={\Kt},{\Ks}$. For $Y={\Kt}$ in \eqref{5.19} we have
\begin{equation*}
\begin{split}
-3&{\eta}_3({\PKt}{\Ks})+{\PoP}^2{\Kt}-{\eta}({\PtP}{\Kt}){\PtK}\\
=&-3{\PKt}+{\PPo}{\PKt}-{\eta}({\Ps}{\Kt}){\phi}{\PsK}\\
=&{\Ks}.
\end{split}
\end{equation*}
Also the both sides hold for $Y={\Ks}$.
\par
Case 3) Put $Y{\in}E_{-1}{\oplus}E_{+1}$. First let us say $Y{\in}E_{-1}$ such that ${\phi}Y=-{\Po}Y$.
The left side of \eqref{5.19} becomes
$$L={\SN}{\ENK}{\PNP}{\phi}^2Y={\PoP}^2Y=-{\Po}Y$$
and the Right side is given by
$$R={\PPo}{\phi}Y={\PoP}{\phi}Y={\Po}{\phi}^2Y=-{\Po}Y.$$
Also the both sides equal to each other for $Y{\in}E_{+1}$.

\vskip 8pt

\begin{remark}\label{Remark 5.2}
\quad In the paper \cite{S2} we have given a complete classification of real hypersurfaces in $G_2({\Bbb
C}^{m+2})$ with commuting Ricci tensor and have proved that they are locally congruent to a tube over a totallly geodesic $\GBo$ in $\GBt$. Also in \cite{S2} we have proved that they have isometric Reeb flows.
\end{remark}
\vskip 6pt

\begin{remark}\label{Remark 5.3}
\quad In the paper \cite{PSW} we have given a complete classification of
pseudo-Einstein hypersurfcaes in $\GBt$ and have found
that there does not exist any Einstein real hypersurface in
$\GBt$.
\end{remark}

\vskip 6pt
\begin{remark}\label{Remark 5.3}
\quad Our main theorem in the introduction will give a contribution to the study of Lie invariant problems for real hypersurfaces in $\NBt$ or to the classification problem of pseudo-Einstein real hypersurfaces in
$\NBt$. Moreover, based on the classification of the isometric Reeb flow in \cite{BS2}, it will give a contribution to the study of real hypersurfaces in complex quadric $Q^m$ with commuting Ricci tensor.  
\end{remark}



\end{document}